\documentclass[11pt]{article}
\usepackage{amssymb,amsfonts,amsmath,latexsym,epsf,tikz,url}
\newtheorem{theorem}{Theorem}[section]

\newtheorem{lemma}[theorem]{Lemma}
\newtheorem{rem}[theorem]{Remark}

\newcommand{\proof}{\noindent{\bf Proof.\ }}
\newcommand{\qed}{\hfill $\square$\medskip}

\begin{document}

\title{Total domination polynomial of graphs from primary subgraphs} 

\author{Saeid Alikhani$^{}$\footnote{Corresponding author} and Nasrin Jafari }

\date{\today}

\maketitle

\begin{center}

   Department of Mathematics, Yazd University, 89195-741, Yazd, Iran\\
{\tt alikhani@yazd.ac.ir}\\

\end{center}

%%%%%%%%%%%%%%ABSTRACT%%%%%%%%%%%%%%%%%%%%%%%%%%%%%%%%%%%%%%%%%%%%%%%%%%%%%%%%%%%%

\begin{abstract}
Let $G = (V, E)$ be a simple graph of order $n$. The total dominating set is a subset $D$ of $V$ that every vertex of $V$ is adjacent to some vertices of $D$. The total domination number of $G$ is equal to minimum cardinality of  total dominating set in $G$ and denoted by $\gamma_t(G)$. The total domination polynomial of $G$ is the polynomial $D_t(G,x)=\sum d_t(G,i)$, where $d_t(G,i)$ is the number of total dominating sets of $G$ of size $i$. Let $G$ be a connected graph constructed from pairwise disjoint connected graphs $G_1,\ldots ,G_k$ by selecting a vertex of $G_1$, a vertex of $G_2$, and identify these two
vertices. Then continue in this manner inductively. We say that $G$ is obtained by point-attaching from $G_1, \ldots ,G_k$ and that $G_i$'s are the primary subgraphs of $G$. 
In this paper, we consider some  particular cases  of these graphs that most of them  are  of importance in chemistry  and study their total domination polynomials.

\end{abstract}


\medskip
\noindent{\bf AMS Subj.\ Class.:} 05A18; 11B 73, 05C12.

%%%%%%%%%%%%%%%%%%%%%%%%%%%%%%%%%%%%%%%%%%%%%%%%%%%%%%%%%%%%%%%%%%%%%%%%%%%%%%%%%
%%%%%%%%%%%%%%%%%%%%%%%%%%%%%%%%%%%%%%%%%%%%%%%%%%%%%%%%%%%%%%%%%%%%%%%%%%%%%%%%%
\section{Introduction}
%%%%%%%%%%%%%%%%%%%%%%%%%%%%%%%%%%%%%%%%%%%%%%%%%%%%%%%%%%%%%%%%%%%%%%%%%%%%%%%%%
%%%%%%%%%%%%%%%%%%%%%%%%%%%%%%%%%%%%%%%%%%%%%%%%%%%%%%%%%%%%%%%%%%%%%%%%%%%%%%%%%

Let $G = (V, E)$ be a simple graph of order $n$. For any vertex $ v \in V$, the open neighborhood of $v$ is the set $N(v)=\{ u \in V | uv \in E\}$ and the closed neighborhood is the set $N[v]=N (v) \cup \{v\}$.
For a set $S\subset V$, the open neighborhood of $S$ is the set $N(S)=\bigcup\limits_{v\in S }N(v)$ and the closed neighborhood of $S$ is the set $N[S]=N (S) \cup S$. A set $D\subset V$ is a total dominating set if every vertex of $V$ is adjacent to some vertices of $D$, or equivalently, $N(D)=V$. The total dominating number $\gamma_t(G)$ is the minimum cardinality of a total dominating set in $G$. A total dominating set with cardinality $\gamma_t(G)$ is called a $\gamma_t$-set. An $i$-subset of $V$ is a subset of $V$ of cardinality $i$. Let $D_t(G, i)$ be the family of total dominating sets of $G$ which are $i$-subsets and let $d_t(G,i)=|D_t(G, i)|$. The polynomial $D_t(G; x)=\sum\limits_{i=1}^n d_t(G,i)x^i$ is defined as total domination polynomial of $G$.

\begin{figure}[ht]%\label{bookroots}
	\hspace{3.5cm}
	\includegraphics[width=7cm]{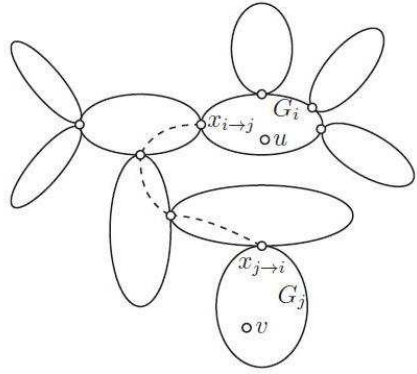}
	\caption{\label{Figure1} Graph $G$ obtained by point-attaching from $G_1,\ldots , G_k$.}
\end{figure}

Let $G$ be a connected graph constructed from pairwise disjoint connected graphs
$G_1,\ldots ,G_k$ as follows. Select a vertex of $G_1$, a vertex of $G_2$, and identify these two vertices. Then continue in this manner inductively.  Note that the graph $G$ constructed in this way has a tree-like structure, the $G_i$'s being its building stones (see Figure \ref{Figure1}).  Usually  say that $G$ is obtained by point-attaching from $G_1,\ldots , G_k$ and that $G_i$'s are the primary subgraphs of $G$. A particular case of this construction is the decomposition of a connected graph into blocks (see \cite{Deutsch}).

\begin{figure}[!h]
 	\begin{center}
 		\includegraphics[width=12cm,height=2cm]{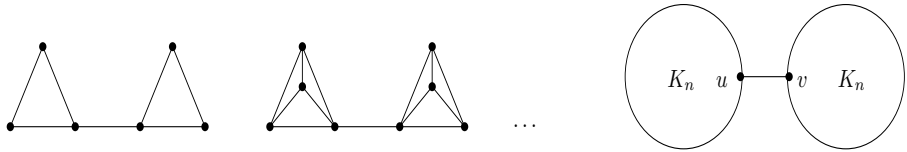}
 		\caption{Barbell graphs $Bar_3, Bar_{4}$ and $Bar_{n}$, respectively.}
 		\label{Bar_n}
 	\end{center}
 \end{figure}

 As an example, the $n$-barbell graph $Bar_n$ with $2n$ vertices, is formed by joining two copies of a complete graph $K_n$ by a single edge (Figure \ref{Bar_n}). Actually, this graph is a specific kind of point-attaching of two complete graphs $K_n$ and the graph $P_2$. Observe that  the total domination polynomial of $n$-barbell graph is 
 $$D_t(Bar_n,x)=\sum\limits_{i=2}^n \binom{2n-2}{i-2}x^i+\sum\limits_{i=n+1}^{2n} \binom{2n}{i}x^i.$$ 
 This formula obtain easily  from counting the total dominating sets of $Bar_n$.     
 
Calculating the total domination polynomial of a graph $G$ is difficult in general, as the smallest power of a non-zero term is the domination number $\gamma_t(G)$ of the graph, and determining whether $\gamma_t(G)\leq  k$ is
known to be NP-complete. So presenting a closed formula for the total domination polynomial of any kind of point-attaching graphs is difficult, but for certain classes of graphs, we can find a closed form expression
for the total domination polynomial.

In this paper, we consider some  particular cases  of these graphs   and study their total domination polynomials.  In Section 2, we consider graphs which obtain by a special point-attaching  of a graph $H$ and $|V(H)|$ copies of graphs $P_3$. We prove that all graphs whose total domination polynomial have just two roots $\{-2,0\}$ are in this form. Also we study the total domination polynomial of some kind of generalized friendship graphs in this section. In Section 3, we investigate the total domination polynomial of cactus chains.

%%% ----------------------------------------------------------------------
\section{Total domination polynomial of graphs from primary subgraphs}

In this section, we consider graphs constructed from primary subgraphs and study their total domination polynomial. Some kind of these graphs have interesting properties. In the Subsection 2.1, we prove that a special kind of graphs from primary subgraphs have exactly two total domination roots. In Subsection 2.2 we study the total domination polynomial of the generalized friendship graph.  

\subsection{Graphs with exactly two total domination roots $\{-2,0\}$}

\begin{figure}
\begin{center}
\includegraphics[width=3.7cm,height=3.2cm]{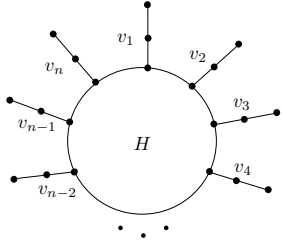}
\caption{The graph $H(3)$.}
\label{H(3)}
\end{center}
\end{figure}

Graphs whose certain  polynomials have
few roots can sometimes give interesting  information about the structure of the graph. 
The characterization of
graphs with few distinct roots of characteristic polynomials
(i.e., graphs with few distinct eigenvalues) have been the
subject of many researchers \cite{bri,dam1,dam2,dam3}. Also the first authors has studied graphs with few domination roots in \cite{few}. 
Let $H$ be an arbitrary graph of order $n$ and consider $n$ copies of graph $P_3$.   By definition, the graph $H(3)$ is obtained by identifying each vertex of $H$ with an end vertex of a $P_3$. See Figure \ref{H(3)}. To obtain the total domination polynomial of $H(3)$, we need the following result. 

\begin{theorem}{\rm\cite{bri}}\label{3}
Let $G$ be a connected graph of order $n\geq 3$. Then $\gamma_ t (G)=\frac{2n}{3}$ if and only if $G$ is $C_3$, $C_6$ or $H(3)$ for some connected graph $H$.
\end{theorem}

\begin{theorem}\label{hop}
For any graph $H$ of order $n$, we have $D_t(H(3),x)=x^{2n}(x+2)^n$.
\end{theorem}
\proof
  Let $D$ be a total dominating set of $H(3)$  of size $k\geq n$ in  Figure \ref{H(3)}. Obviously $\{v_1,v_2,\ldots,v_n\}\subset D$. To choose $n+i$ ($0\leq i\leq n$) other vertices of $V(H(3))\setminus \{v_1,v_2,\ldots,v_n\}$, we have $\binom{n}{i}2^{n-i}$ possibilities. So we have the result.\qed
  
  Now, we state and  prove the following result.  
  \begin{theorem}\label{2roots}
  Let $G$ be a graph. Then $D_t(G, x)=x^{2n}(x+2)^n$ if and only if $G= H(3)$ for some graph $H$ of order $n$.
  \end{theorem}
  \proof
   ($\Leftarrow$) It follows from  Theorem \ref{hop}.
   
   ($\Rightarrow$) Let $G$ be a graph with $D_t(G, x)=x^{2n}(x+2)^n$. Thus $|V(G)|=3n$ and $G$ has no isolated vertex. Since $\gamma_t(G)=2n$, by Theorem \ref{3}, every  component of   $G$ is a cycle $C_3$, $C_6$ or a $H(3)$ for some connected graph $H$. Since $D_t(C_3,x)=x^2(x+3)$   and $D_t(C_6,x)=x^4(x+3)^2$ does not divide $x^{2n}(x+1)^n$, we conclude that there exists a graph $H$ such that $G=H(3)$ and the proof is complete. \qed

\begin{rem}   
      The characterization of graphs whose graph polynomials have few roots have been an interesting problem and studied well in the literature (\cite{few,Ebrahim}). Also there is a conjecture
  in  \cite{nasrin} which states that every integer total domination roots is in the set $\{-3,-2,-1,0\}$. So finding the graphs whose total domination polynomial  have these few roots can be a good start for solving this conjecture. Theorem \ref{2roots} characterize all graphs whose total domination polynomial have just two distinct roots $-2$ and $0$.    
  \end{rem}  
   
   \medskip
   
   \subsection{Total domination polynomial of the generalized friendship graph}
\begin{figure}
\begin{center}
\includegraphics[width=12cm,height=3cm]{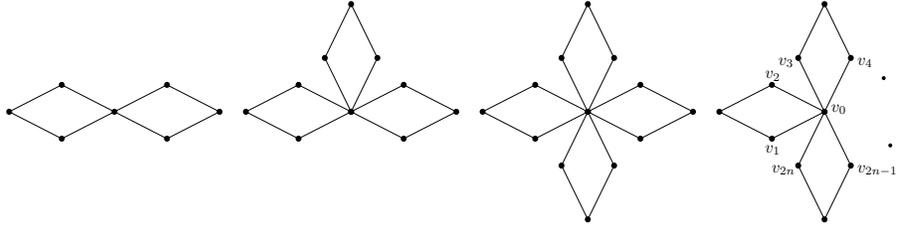}
\caption{Friendship graphs $F_{1,4}, F_{2,4}, F_{3,4}$ and $F_{n,4}$, respectively.}
\label{Dutch}
\end{center}
\end{figure}

  Now we consider another kind of point-attaching graphs and study their total domination polynomials.  
  The generalized friendship graph $F_{n,q}$ is a collection of $n$ cycles (all of order $q$), meeting at a common vertex (see Figure \ref{Dutch}). The generalized friendship graph may also be referred to as a flower \cite{Schiermeyer}. For $q=3$ the graph $F_{n,q}$ is denoted simply by $F_n$ and is friendship graph. The total domination polynomial of $F_n$ and its roots studied in \cite{nasrin}. Here, first we compute the total domination number of $F_{n,4}$. 

\begin{theorem}
	 For any $n\geq 1$, we have  $\gamma_t(F_{n,4})=n+1$.
\end{theorem}
\proof
  Let $\{v_1,\ldots ,v_{2n}\}$ be vertex set of $F_{n,4}$ that adjacent by $v_0$ (common vertex in $F_{n,4}$). Then $\{v_0, v_1,v_3,\ldots , v_{2n-1}\}$ is a total dominating set for $F_{n,4}$ (see Figure \ref{Dutch}) and the set  $D\subseteq V(F_{n,4})$ of size less than or equal  $n$ is not total dominating set for $F_{n,4}$ , therefore $\gamma_t(F_{n,4})=n+1$.\qed
  
The following theorem is useful for finding the recurrence relations of the domination polynomials of graphs. The vertex contraction $G/u$ of a graph $G$ by a vertex $u$ is the operation under which all vertices in $N(u)$
are joined to each other and then $u$ is deleted (see \cite{Walsh}). 

\begin{theorem}{\rm\cite{Dod}}  \label{thm}
           \begin{enumerate}
      \item[(i)]For any vertex $u$ in the graph $G$ we have 
      \begin{center}
      $D_t(G, x)=D_t(G\setminus u, x)+xD_t(G/u, x)+x^2\sum \limits _{v\in N(u)}D_t(G\setminus N[\{u,v\}],x)-(1+x)p_u(G)$,
      \end{center}
       where $p_u(G,x)$ is the polynomial counting the total dominating sets of $G\setminus u$ which do not contain any vertex of $N(u)$ in $G$.
      \item[(ii)]Let $u,v \in V(G)$ be two non-adjacent vertices of $G$ with $N(v)\subseteq N(u)$. Then
        \begin{center}
      $D_t(G, x)=D_t(G\setminus u, x)+xD_t(G/u, x)+x^2\sum \limits _{w\in N(u)\cap N(v)}D_t(G\setminus N[\{u,w\}],x)$.
      \end{center}
      \item[(iii)]If $u,v \in V(G)$ be two vertices of $G$ with $N[v]\subseteq N[u]$. Then
        \begin{center}
      $D_t(G, x)=D_t(G\setminus u, x)+xD_t(G/u, x)+x^2\sum \limits _{w\in N(u)}D_t(G\setminus N[\{u,w\}],x)$.
      \end{center}
      \item[(iv)]Let  $e=\{u,v\} \in E(G)$ and $N[v]=N[u]$. Then
        \begin{center}
      $D_t(G, x)=D_t(G\setminus e, x)+x^2D_t(G\setminus N[u],x)$.
      \end{center}
      \end{enumerate}
\end{theorem}
  
  Now we state and prove a recurrence relation for the total domination polynomial of $F_{n,4}$. 
\begin{theorem}\label{Fn4}
 For any $n \geq 1$, we have
 \begin{center}
 $D_t(F_{n,4},x)=x(x+2)[(x+1)D_t(F_{n-1,4},x)-(x^3+2x)^{n-1}]$,
\end{center}
with initial value $D_t(F_{1,4},x)=x^4+4x^3+4x^2$.
\end{theorem}
\proof
	Consider graph $F_{n,4}$ and $u$, $v$ as shown in Figure \ref{f_n,4}. By Theorem \ref{thm}, we have 
        \begin{align*}
               &D_t(F_{n,4},x)\stackrel{part (ii)}{=}D_t(F_{n,4}\setminus u,x)+xD_t(F_{n,4}/u,x)+x^2D_t(P_3,x)^{n-1}
               \\&\hspace*{1.8cm}\stackrel{part (iii)}{=}\overbrace{D_t((F_{n,4}\setminus u)\setminus v,x)}^{0}+xD_t((F_{n,4}\setminus u)/v,x)+x^2D_t(P_3,x)^{n-1}
               \\&\hspace*{2.5cm}+x[D_t((F_{n,4}/u)\setminus v,x) +xD_t((F_{n,4}/u)/v,x)+x^2D_t(P_3,x)^{n-1}]
                \\&\hspace*{2.5cm}+x^2D_t(P_3,x)^{n-1}
               \\&\hspace*{2.2cm}=(x^2+2x)D_t((F_{n,4}\setminus u)/v,x)+(x^3+2x^2)D_t(P_3,x)^{n-1}
                \\&\hspace*{1.9cm}\stackrel{part (i)}{=}(x^2+2x)[D_t(F_{n-1,4},x)+xD_t(F_{n-1,4},x)-(x+1)P_w((F_{n,4}\setminus u)/v,x)]
                \\&\hspace*{2.5cm}+(x^3+2x^2)D_t(P_3,x)^{n-1}
               \\&\hspace*{2.1cm}=x(x+1)(x+2)D_t(F_{n-1,4},x)-x(x+1)(x+2)P_w((F_{n,4}\setminus u)/v,x)
               \\&\hspace*{2.5cm}+(x^3+2x^2)D_t(P_3,x)^{n-1}.
        \end{align*}
Since $P_w((F_{n,4}\setminus u)/v,x)=D_t(P_3,x)^{n-1}$, so we have the result.\qed

\begin{figure}
\begin{center}
\includegraphics[width=13cm,height=4.3cm]{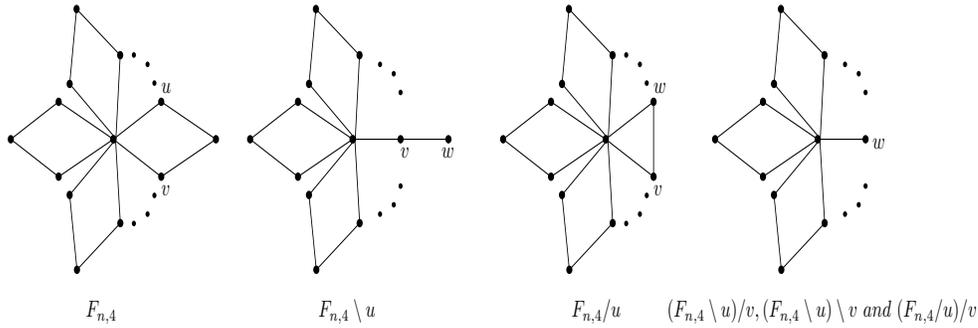}
\caption{Graphs which considered in Theorem \ref{Fn4}}
\label{f_n,4}
\end{center}
\end{figure}

\section{Total domination polynomial of cactus chains}
In this section, we consider another kind of point-attaching graphs and study their total domination
polynomials. These kind of graphs are important in Chemistry. A cactus graph is a connected graph in which no edge lies in more than one cycle. Consequently, each block of a cactus graph is either an edge or a cycle. If all blocks of a cactus $G$ are cycles of the same size $i$ , the cactus is $i$ -uniform. A triangular cactus is a graph whose blocks are triangles, i.e., a $3$ -uniform cactus. We call the number of triangles in $G$ the length of the chain. Obviously, all chain triangular cactus of the same length are isomorphic. Hence, we denote the chain triangular cactus of length $n$ by $T_n$ (see Figure \ref{Gn}). By replacing triangles in this definitions by cycles of length $4$ we obtain cacti whose every block is $C_4$. In this section we shall study the total domination polynomial of cactus chains. 

\subsection{Total domination polynomial of the chain triangular cactus}
In this subsection we shall study the total domination polynomial of chain triangular cactus. To do this, we consider graph $G_n$ as shown in Figure \ref{Gn}, which is also a kind of point attaching graphs. Note that $G_n$ is a point attaching of $T_n$ and $P_2$.  First we state and prove the following theorem:  

%\begin{figure}[!h]
%\begin{center}
%\includegraphics[width=4.5cm,height=1.3cm]{T_n}
%\caption{The chain triangular cactus}
%\label{T_n}
%\end{center}
%\end{figure}

\begin{figure}[!h]
\begin{center}
\includegraphics[width=11cm,height=1.5cm]{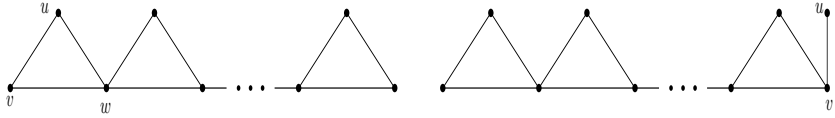}
\caption{The graphs $T_n$ and $G_n$, respectively.}
\label{Gn}
\end{center}
\end{figure}

\begin{theorem}
	For every $n\geq 2$, 
	\begin{center}
	   $D_t(G_n,x)=(x+1)[D_t(T_n,x)-D_t(G_{n-1},x)]+x^2D_t(G_{n-2},x),$
	   \end{center}
	   where $D_t(G_0,x)=x^2$, $D_t(G_1,x)=x^4+3x^3+3x^2$ and $D_t(T_2,x)=x^5+5x^4++6x^3+4x^2$.
\end{theorem}
\proof
	Consider the graph $G_n$ as shown in Figure \ref{Gn}. Since $G_n\setminus u$ is isomorphic to $G_n/u$, by Theorem \ref{thm}(i), we have
	\begin{align*}
	&D_t(G_n,x)=(x+1)D_t(G_n/u,x)+x^2D_t(G_{n}\setminus N[\{u,v\}],x)-(x+1)p_u(G_n)
	\\&\hspace*{1.6cm}=(x+1)D_t(T_n,x)+x^2D_t(G_{n-2},x)-(x+1)D_t(G_{n-1},x).
	\end{align*}
  Note that $p_u(G_n)=D_t(G_{n-1},x)$.
	\qed

\begin{theorem}
  For every $n\geq 3$, 
	 \begin{center}
	   $D_t(T_n,x)=(x+1)D_t(G_{n-1},x)+x^2[D_t(G_{n-2},x)+D_t(G_{n-3},x)]$.
	   \end{center}
\end{theorem}   
\proof
     Consider the graph $T_n$ and its vertex $u$ as shown in the Figure \ref{Gn}. By Theorem \ref{thm}(iii), we have
	\begin{align*}
	&D_t(T_n,x)=(x+1)D_t(T_n/u,x)+x^2[D_t(T_{n}\setminus N[\{u,v\}],x)+D_t(T_{n}\setminus N[\{u,w\}],x)]
	\\&\hspace*{1.5cm}=(x+1)D_t(G_{n-1},x)+x^2[D_t(G_{n-2},x)+D_t(G_{n-3},x)].
	\end{align*}
     Note that $T_n\setminus u$ is isomorphic to $T_n/u$.\qed

    \subsection{Total domination polynomial of para-chain square cactus graphs}
    
\begin{figure}[!h]
    \begin{center}
    \includegraphics[width=10cm,height=1.7cm]{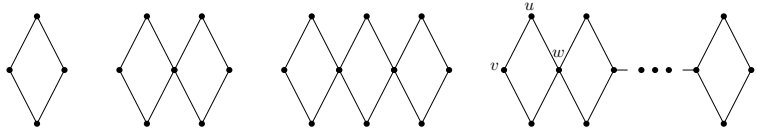}
    \caption{Para-chain square cactus graph $Q_{1}, Q_{2}, Q_{3}$ and $Q_{n}$, respectively.}
    \label{Q_n}
    \end{center}
    \end{figure}

    In this subsection we consider a para-chain of length $n$ , $Q_n$ , as shown in Figure \ref{Q_n} and obtain a recurrence
relation for the domination polynomial of $Q_n$.
    
    \begin{figure}[!h]
    \begin{center}
    \includegraphics[width=10cm,height=3.4cm]{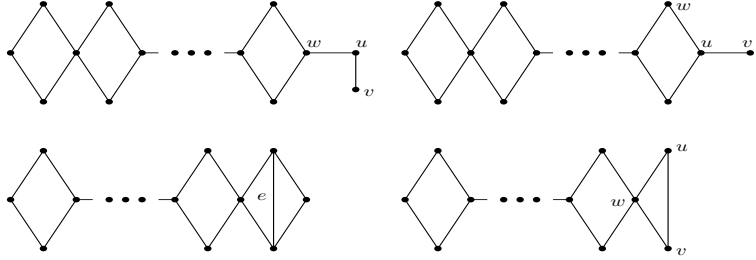}
    \caption{graphs $Q_{n}(2), Q_n(1), Q_n+e$ and $Q_{n}^\Delta$, respectively.}
    \label{Qn}
    \end{center}
    \end{figure}
    
    \begin{figure}[!h]
    \begin{center}
    \includegraphics[width=5cm,height=1.7cm]{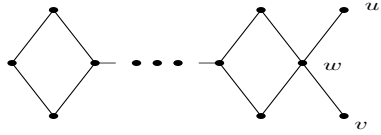}
    \caption{graph $Q'_{n}$.}
    \label{Q'n}
    \end{center}
    \end{figure}
    
\begin{lemma}\label{le}
    For graphs in Figures \ref{Qn} and \ref{Q'n} have:
    \begin{enumerate}
    \item[(i)]$D_t(Q_n(1),x)=xD_t(Q_{n},x)+x^2[D_t(Q_{n-1},x)+3D_t(Q^{'}_{n-1},x)]$, where $D_t(Q_0(1),x)=x^2$.
    \item[(ii)]$D_t(Q_n(2),x)=x^2[D_t(Q_n,x)+(x+1)D_t(Q_{n-1},x)+(3x+1)D_t(Q^{'}_{n-1},x)]$, where $D_t(Q_0(2),x)=x^3+2x^2$.
    \item[(iii)]$D_t(Q^{'}_n,x)=x(x+1)D_t(Q_{n},x)+x^2(x+2)D_t(Q_{n-1},x)+3x^2(x+1)D_t(Q^{'}_{n-1},x)$, where $D_t(Q_0^{'},x)=x^3+2x^2$.
    \item[(iv)]$D_t(Q_n^\Delta,x)=x(x+1)D_t(Q_n,x)+x^2(x+2)D_t(Q_{n-2},x)+x^2(3x+4)D_t(Q^{'}_{n-1},x)$, where $D_t(Q_0^\Delta,x)=x^3+3x^2$.
    \end{enumerate}
\end{lemma}
\proof
	\begin{enumerate}
	\item[(i)]
 Consider graph $Q_n(1)$ and its vertex $u$  in Figure \ref{Qn}. By Theorem \ref{thm}(iii),
          \begin{align*}
          & D_t(Q_n(1),x)=\overbrace{D_t(Q_n(1)\setminus u,x)}^{0}+x D_t(Q_n(1)/u,x)
	    \\&\hspace*{2.4cm}+x^2[D_t(Q_n(1)\setminus N[\{u,v\}],x)+2D_t(Q_n(1)\setminus N[\{u,w\}],x)]
	    \\& \hspace*{2.1cm}=xD_t(Q_{n}+e,x)+x^2[D_t(Q_{n-1}x)+2D_t(Q^{'}_{n-1},x)].
        \end{align*}
         By applying Theorem \ref{thm}(iv) on $Q_n+e$, we have
           \begin{center}
           $D_t(Q_n+e,x)=D_t(Q_{n},x)+x^2D_t(Q^{'}_{n-1},x)$.
           \end{center}
         So we have result.
	  
      \item[(ii)]
      Consider the vertex  $u$  as  shown in Figure \ref{Qn}. By Theorem \ref{thm}(iii), we have
	   \begin{align*}
	    &D_t(Q_n(2),x)=\overbrace{D_t(Q_n(2)\setminus u,x)}^{0}+x D_t(Q_n(2)/u,x)
	    \\&\hspace*{2.4cm}+x^2[D_t(Q_n(2)\setminus N[\{u,v\}],x)+D_t(Q_n(2)\setminus N[\{u,w\}],x)]
	    \\& \hspace*{2.1cm}=xD_t(Q_{n-1}(1),x)+x^2[D_t(Q^{'}_{n-1},x)+D_t(Q_{n-1},x)].
	  \end{align*}
        By using Part $(i)$ in the above equation, we have result.
      \item[(iii)]
        Consider graph $Q^{'}_n$ in Figure \ref{Q'n}. By Theorem \ref{thm}, 
         \begin{align*}
         & D_t(Q^{'}_n,x)=D_t(Q^{'}_n\setminus u,x)+x D_t(Q^{'}_n/u,x)+x^2D_t(Q^{'}_n\setminus N[\{u,w\}],x)
	    \\& \hspace*{1.7cm}=(x+1)D_t(Q_{n}(1),x)+x^2D_t(Q_{n-1},x).
         \end{align*}
        Using Part $(i)$ in the above equation, we have result.
     \item[(iv)]
       By Theorem \ref{thm}(iii), we have
          \begin{align*}
          & D_t(Q_n^\Delta ,x)=D_t(Q_n^\Delta \setminus u,x)+x D_t(Q_n^\Delta /u,x)
	    \\&\hspace*{2.1cm}+x^2[D_t(Q_n^\Delta \setminus N[\{u,v\}],x)+D_t(Q_n^\Delta \setminus N[\{u,w\}],x)]
	    \\& \hspace*{1.7cm}=(x+1)D_t(Q_{n}(1),x)+x^2[D_t(Q^{'}_{n-1},x)+D_t(Q_{n-1},x)].
          \end{align*}
       So by using Part $(i)$ in the above equation, we have result.\qed
   \end{enumerate}
   
 \begin{theorem}
 The total domination polynomial of para-chain $Q_n$ is given by
     \begin{align*}
      &D_t(Q_n,x)=x^2(x+2)[D_t(Q_{n-1},x)+D_t(Q_{n-2},x)+xD_t(Q_{n-3},x)]
      \\&\hspace*{2cm}+x^2(3x^2+7x+2)D_t(Q^{'}_{n-2},x),
    \end{align*}
    where $D_t(Q_1,x)=x^4+4x^3+4x^2$ and $D_t(Q_2,x)=D_t(F_{2,4},x)$.
 \end{theorem}
    \proof
        With regards to Figure \ref{Q_n} and Theorem \ref{thm}, we have 
        \begin{align*}
        &D_t(Q_n,x)=D_t(Q_n\setminus u,x)+xD_t(Q_n/u,x)
      \\&\hspace*{2cm}+x^2[D_t(Q_n\setminus N[\{u,v\}],x)+D_t(Q_n\setminus N[\{u,w\}],x)
      \\&\hspace*{1.6cm} =D_t(Q_{n-1}(2),x)+xD_t(Q_{n-1}^\Delta ,x)+x^2[D_t(Q_{n-2}^{'},x)+D_t(Q_{n-2},x)]
        \end{align*}
    Now by Lemma \ref{le} results is obtained.\qed

 \subsection{Total domination polynomial of ortho-chain square cactus graphs}
   
\begin{figure}[!h]
    \begin{center}
    \includegraphics[width=10cm,height=1.7cm]{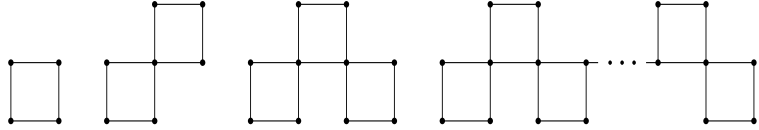}
    \caption{ortho-chain square cactus graphs $O_{1}, O_{2}, O_{3}$ and $O_{n}$, respectively.}
    \label{O_n}
    \end{center}
    \end{figure}

   In this subsection we consider an ortho-chain of length $n$, $O_n$, as shown in Figure \ref{O_n}. We shall obtain a
recurrence relation for the total  domination polynomial of $O_n$.

    \begin{figure}[!h]
    \begin{center}
    \includegraphics[width=10cm,height=3.4cm]{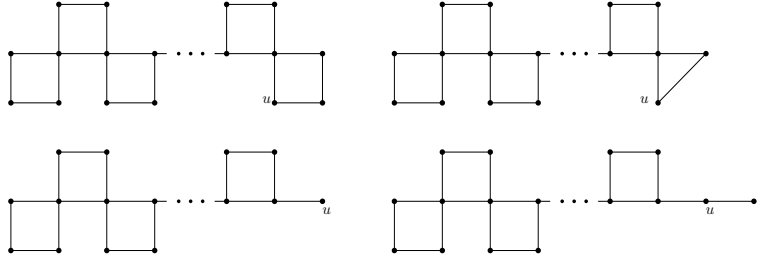}
    \caption{graphs $O_n$, $O_{n}^\Delta$, $O_n(1)$ and $O_{n}(2)$, respectively.}
    \label{On}
    \end{center}
    \end{figure}

  \begin{lemma}\label{lee}
   For graphs in Figure \ref{On}, we have
   \begin{enumerate}
    \item[(i)]$D_t(O_n(1),x)=(x+1)[D_t(O_{n},x)-D_t(O_{n-1}(2),x)]$, where $D_t(O_0(1),x)=x^2$.
    \item[(ii)]$D_t(O_n(2),x)=x[(x+1)D_t(O_{n},x)-D_t(O_{n-1}(2),x)]$, where $D_t(O_0(2),x)=x^3+2x^2$.
    \item[(iii)]$D_t(O_n^\Delta,x)=(x+1)^2D_t(O_n,x)-(2x+1)D_t(O_{n-1}(2),x)$, where $D_t(O_0^\Delta ,x)=x^3+3x^2$.
   \end{enumerate}
  \end{lemma}
   
   \proof
   % \begin{enumerate}
         Similar to lemma \ref{le} by using theorem \ref{thm} we have result. Note that in part $(i)$, $P_u(O_n)=D_t(O_{n-1}(2),x)$.
%    	\item[(i)]
%        The proof follows from Theorem \ref{thm} for vertex $u$ in graph $O_n(1)$. 
%        
%        \item[(ii)]         The proof follows from Theorem \ref{thm} for vertex $u$ in graph $O_n(2)$.
%        
%         \item[(iii)]         The proof follows from Theorem \ref{thm} for vertex $u$ in graph and $O_n^\Delta$ in Figure \ref{On}. \qed
%       \end{enumerate}  
        
        \begin{theorem}
          For graph $O_n$ in Figure \ref{O_n}, we have
           \begin{center}
            $D_t(O_n,x)=x(x+2)[(x+1)D_t(O_{n-1},x)-D_t(O_{n-2}(2),x)]$,
           \end{center}
        where $D_t(O_1,x)=x^4+4x^3+4x^2$ and $D_t(O_2,x)=D_t(F_{2,4},x)$.
        \end{theorem}
        
        \proof
         Consider graph $O_n$ and vertex $u$ as shown in Figure \ref{On}. By Theorem \ref{thm} for vertex $u$, we have 
         \begin{center}
         $D_t(O_n,x)=D_t(O_{n-1}(2),x)+xD_t(O_{n-1}^\Delta ,x)+x^2D_t(O_{n-2}(2),x)$.
         \end{center}
         Now by using Lemma \ref{lee}(ii) and (iii), we have result.\qed

%----------------------------------------------------------------------------

% ------------------------------------------------------------------------

% ------------------------------------------------------------------------
\end{document}